\documentclass{conm-p-l}

\usepackage{amsfonts,amsmath,amssymb,verbatim,graphics,axodraw}

\theoremstyle{plain}
\newtheorem{theorem}{Theorem}[section]

\theoremstyle{definition}
\newtheorem{definition}[theorem]{Definition}
\newtheorem{example}[theorem]{Example}
\newtheorem{conjecture}[theorem]{Conjecture}

\numberwithin{equation}{section}

\newcommand{\qbin}[2]{\genfrac{[}{]}{0pt}{}{#1}{#2}}
\newcommand{\qbins}[2]{{\textstyle\genfrac{[}{]}{0pt}{}{#1}{#2}}}
\newcommand{\Z}{\mathbb{Z}}

\newcommand{\N}{\mathbb{N}}

\newcommand{\CC}{\mathbb{C}}
\newcommand{\R}{\mathbb{R}}

\newcommand{\gggg}{\mathfrak{g}}
\newcommand{\gggh}{\hat{\gggg}}
\newcommand{\wt}{\mathrm{wt}}
\newcommand{\lev}{\mathrm{lev}}

\newcommand{\la}{\lambda}
\newcommand{\lt}{\tilde{\ell}}
\newcommand{\La}{\Lambda}

\renewcommand{\P}{\mathcal{P}}
\newcommand{\Pbar}{\overline{\mathcal{P}}}
\newcommand{\D}{\mathcal{D}}
\newcommand{\Ebar}{\overline{E}}
\newcommand{\Dbar}{\overline{D}}

\newcommand{\Sbar}{\overline{S}}
\newcommand{\Xbar}{\overline{X}}
\newcommand{\Fbar}{\overline{F}}
\newcommand{\Wh}{\hat{W}}

\newcommand{\sg}{\sigma_{\mathrm{G}}}
\newcommand{\ch}{c}
\newcommand{\cc}{cc}
\newcommand{\inner}[2]{\langle #1\,,\,#2\rangle}

\newcommand{\RC}{\mathrm{RC}}
\newcommand{\CST}{\mathrm{CST}}
\newcommand{\SCST}{\mathrm{SCST}}

\newcommand{\Image}{\mathrm{Im}}

\begin{document}

\title{Crystal bases and $q$-identities}

\author[M.~Okado]{Masato Okado}
\address{Department of Informatics and Mathematical Science, 
Graduate School of Engineering Science, Osaka University, 
Toyonaka, Osaka 560-8531, Japan}
\email{okado@sigmath.es.osaka-u.ac.jp}

\author[A.~Schilling]{Anne Schilling}
\address{Department of Mathematics 2-279, M.I.T., Cambridge, MA 02139, U.S.A.
and\newline
Department of Mathematics, University of California, One Shields 
Avenue, Davis, CA 95616-8633, U.S.A.}
\email{anne@math.mit.edu, anne@math.ucdavis.edu}

\author[M.~Shimozono]{Mark Shimozono}
\address{Department of Mathematics, 460 McBryde Hall, Virginia Tech,
Blacksburg, VA 24061-0123, U.S.A}
\email{mshimo@math.vt.edu}



\begin{abstract}
The relation of crystal bases with $q$-identities is discussed,
and some new results on crystals and $q$-identities associated
with the affine Lie algebra $C_n^{(1)}$ are presented.
\end{abstract}

\maketitle

\tableofcontents

\section{Introduction}

The purpose of this paper is two-fold. First,
we would like to advocate the importance
of crystal theory to the theory of $q$-series.
In particular crystal base theory provides a unifying and general 
setting for a large class of $q$-identities.
Second, as evidence, some new identities associated to the affine Lie algebra 
$C_n^{(1)}$ are presented. The emphasis here will not be on
the completeness of the results since the field is evolving quite rapidly, 
but rather on the presentation of the main ideas and techniques
used.

The Rogers--Ramanujan identities are undoubtably the most
famous $q$-series identities. They are given by
\begin{align}
\label{RR1}
\sum_{n=0}^\infty \frac{q^{n^2}}{(q)_n} &=
 \prod_{j=1}^\infty \frac{1}{(1-q^{5j-4})(1-q^{5j-1})}\\
\label{RR2}
\sum_{n=0}^\infty \frac{q^{n(n+1)}}{(q)_n} &=
 \prod_{j=1}^\infty \frac{1}{(1-q^{5j-3})(1-q^{5j-2})}
\end{align}
where $(q)_n=(1-q)(1-q^2)\cdots(1-q^n)$. We will view them
here as identities of formal power series, meaning that
in the expansion as series in the formal variable $q$
the coefficients of $q^N$ match on both sides for all $N\ge 0$.
What contributes to their beauty is that these coefficients
can be interpreted combinatorially. The coefficient of $q^N$
on the left-hand side of \eqref{RR1} is the number of partitions
of $N$ for which the difference between any two parts is at least
two. The coefficient of $q^N$ of the right-hand side
of \eqref{RR1} on the other hand is the number of partitions of $N$
with parts congruent to 1 or 4 modulo 5.
Similarly, the coefficients of $q^N$ on the left and right side 
of \eqref{RR2} can be interpreted as the number of partitions
of $N$ for which the difference between any two parts is at least
two and the smallest part is greater than one, and the number of
partitions of $N$ with parts congruent to 2 and 3 modulo 5, respectively.

Many of the ideas regarding crystals and $q$-identities can already
be demonstrated in terms of the Rogers--Ramanujan identities.
The point of focus here shall be the debut of the Rogers--Ramanujan
identities on the mathematical physics stage, in particular their
appearance in the hard hexagon model in a paper by Baxter \cite{B81} in 1981.
In this setting the Rogers--Ramanujan identities can be viewed
as two different evaluations of the generating function of certain paths
which are coined bosonic and fermionic evaluations.
Details are discussed in section \ref{sec:HHM}.
As it turns out the relation between the Rogers--Ramanujan identities
and the hard hexagon model is only part of a much bigger picture.
In terms of representation theory, the paths that occur in the hard 
hexagon model are elements of tensor products of crystals
associated with the affine Lie algebra $\hat{\mathfrak{sl}}_2$.
Crystal bases were introduced by Kashiwara \cite{K1,K2} and 
roughly speaking are bases of representations of 
quantum universal enveloping algebras $U_q(\gggg)$ as the parameter
$q$ (not to be confused with the $q$ in the $q$-series!) tends to zero.
Here $\gggg$ is any symmetrizable Kac--Moody algebra.
As in the Rogers--Ramanujan case, there are two different ways
to evaluate generating functions of tensor products of crystals,
thereby giving rise to $q$-identities.
Hence crystal base theory provides a natural framework for $q$-identities.
Crystal bases, path spaces and their generating functions
are discussed in section \ref{sec:crystal}. The two different
ways to evaluate the paths generating functions are subject of
sections \ref{sec:bose} and \ref{sec:fermi}, respectively.
In particular in section \ref{sec:lev C} we present new fermionic
formulas for level-restricted paths.
We close in section \ref{sec:open} with some outstanding open problems.

Before indulging in the fascinating theory of crystal
bases there is one important point that needs to be addressed.
Applying Jacobi's triple product identity
\begin{equation*}
\sum_{n=-\infty}^\infty z^n q^{n^2} =
 \prod_{n=0}^\infty (1-q^{2n+2})(1+z q^{2n+1})(1+z^{-1} q^{2n+1})
\end{equation*}
the right-hand sides of \eqref{RR1} and \eqref{RR2} can be rewritten as
alternating sums yielding the identities
\begin{align}
\label{bf1}
\sum_{n=0}^\infty \frac{q^{n^2}}{(q)_n} &=
 \frac{1}{(q)_\infty} \sum_{j=-\infty}^\infty (-1)^j q^{\frac{j}{2}(5j+1)}\\
\label{bf2}
\sum_{n=0}^\infty \frac{q^{n(n+1)}}{(q)_n} &=
 \frac{1}{(q)_\infty} \sum_{j=-\infty}^\infty (-1)^j q^{\frac{j}{2}(5j+3)}.
\end{align}
The two different evaluations of generating functions of crystal paths
that were mentioned above really yield (polynomial) analogues of \eqref{bf1} 
and \eqref{bf2} rather than \eqref{RR1} and \eqref{RR2}.
Identities relating sums to alternating sums are more general
than identities relating sums to products. Only in special
cases can the alternating sums be evaluated as products,
namely when Jacobi's triple product identity or more generally
the Macdonald identities \cite{Mac} can be applied.

\subsection*{Acknowledgements}
We are grateful to Tim Baker for discussions.
A.S. would like to thank the organizers Bruce Berndt and Ken Ono 
for the invitation to present these results
at the conference on $q$-Series with Applications to Combinatorics, 
Number Theory, and Physics held in Urbana-Champaign in October 2000.
Thanks for this exciting conference!

\section{The Rogers--Ramanujan identities and the Hard Hexagon model}
\label{sec:HHM}

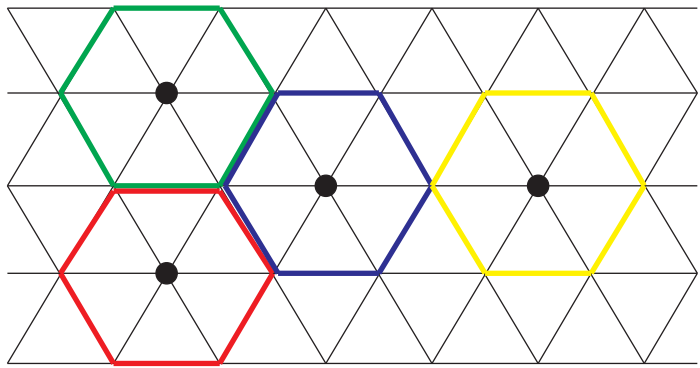
\begin{figure}
\begin{picture}(260,140)(0,0)
\Line(0,0)(260,0)
\Line(0,34)(260,34)
\Line(0,67)(260,67)
\Line(0,102)(260,102)
\Line(0,134)(260,134)
\Line(0,67)(40,134)
\Line(0,0)(80,134)
\Line(40,0)(120,134)
\Line(80,0)(160,134)
\Line(120,0)(200,134)
\Line(160,0)(240,134)
\Line(200,0)(260,102)
\Line(240,0)(260,34)
\Line(40,0)(0,67)
\Line(80,0)(0,134)
\Line(120,0)(40,134)
\Line(160,0)(80,134)
\Line(200,0)(120,134)
\Line(240,0)(160,134)
\Line(260,34)(200,134)
\Line(260,102)(240,134)
\CCirc(60,102){4}{Black}{Black}
\CCirc(60,34){4}{Black}{Black}
\CCirc(120,67){4}{Black}{Black}
\CCirc(200,67){4}{Black}{Black}
\SetWidth{2}
\SetColor{Green}
\Line(40,67)(80,67)
\Line(80,67)(100,102)
\Line(100,102)(80,134)
\Line(80,134)(40,134)
\Line(40,134)(20,102)
\Line(20,102)(40,67)
\SetColor{Red}
\Line(40,0)(80,0)
\Line(80,0)(100,34)
\Line(100,34)(80,65)
\Line(80,65)(40,65)
\Line(40,65)(20,34)
\Line(20,34)(40,0)
\SetColor{Blue}
\Line(102,34)(140,34)
\Line(140,34)(160,67)
\Line(160,67)(140,102)
\Line(140,102)(102,102)
\Line(102,102)(82,67)
\Line(82,67)(102,34)
\SetColor{Yellow}
\Line(180,34)(220,34)
\Line(220,34)(240,67)
\Line(240,67)(220,102)
\Line(220,102)(180,102)
\Line(180,102)(160,67)
\Line(160,67)(180,34)
\end{picture}
\caption{An allowed configuration of particles in the hard hexagon model.
The hexagons around the particles are drawn in different colours.}
\label{fig:hhconfig}
\end{figure}

The hard hexagon model is a two-dimensional lattice model of a gas
of hard or non-overlapping particles. The particles are placed
on a triangular lattice such that no two particles can occupy
two adjacent sites. If one views each particle as the center of 
a hexagon, then the condition that no two particles can be adjacent
translates into the condition that the hexagons cannot overlap.
This explains the name of the model.
An example of an allowed particle configuration is
shown in Figure \ref{fig:hhconfig}.

To pack the lattice densely all particles must lie on 
one of the three sublattices corresponding to the three corners of 
the triangles of the lattice. Hence one of the sublattices
is distinguished from the others. At low particle density the probability
for a particle to be on a particular site is equal for sites
on all three sublattices (assuming either an infinite or sufficiently
large lattice so that boundary effects can be neglected). 
Let $\rho_a$ for $a=1,2,3$ be the probability that there is a particle 
at a fixed site on sublattice $1,2,3$, respectively.
If the boundary conditions of the model are such that at close packing
all particles are on sublattice 1 then it is intuitively clear that
the order parameter defined as $R=\rho_1-\rho_2$ must undergo a phase 
transition. At low densities $R$ is zero, but at some critical
density $R$ will become positive until at high densities it is one.
Baxter \cite{B81} managed to determine the precise point at which
the phase transition occurs exactly by using corner transfer matrices.
In essence, the corner transfer matrix method reduces the
two-dimensional problem to a one-dimensional problem. The precise
details are beyond the scope of this paper and can be found in
section 14 of Baxter's book \cite{B82}.

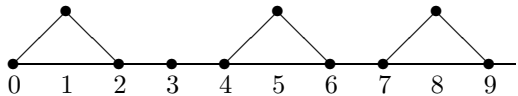
\begin{figure}
\begin{picture}(200,50)(0,0)
\put(10,20){\line(1,1){20}}
\put(50,20){\line(-1,1){20}}
\put(90,20){\line(1,1){20}}
\put(130,20){\line(-1,1){20}}
\put(150,20){\line(1,1){20}}
\put(190,20){\line(-1,1){20}}
\put(8,9){0}
\put(28,9){1}
\put(48,9){2}
\put(68,9){3}
\put(88,9){4}
\put(108,9){5}
\put(128,9){6}
\put(148,9){7}
\put(168,9){8}
\put(188,9){9}
\put(10,20){\line(1,0){190}}
\put(10,20){\circle*{4}}
\put(30,40){\circle*{4}}
\put(50,20){\circle*{4}}
\put(70,20){\circle*{4}}
\put(90,20){\circle*{4}}
\put(110,40){\circle*{4}}
\put(130,20){\circle*{4}}
\put(150,20){\circle*{4}}
\put(170,40){\circle*{4}}
\put(190,20){\circle*{4}}
\end{picture}
\caption{A path of length $9$}
\label{fig path}
\end{figure}

The one-dimensional problem that Baxter encountered and which turns
out to be of importance to the Rogers--Ramanujan identities
is the following.
Consider $L+1$ points on a line labeled by $i=0,1,2,\ldots,L$.
Assign to each point a height variable $\sigma_i$ which
takes on the values $0$ or $1$. In addition the height variables 
satisfy the restrictions $\sigma_0=\sigma_L=0$ and 
$\sigma_i\sigma_{i+1}=0$. An allowed configuration of height variables
for a given length $L$ is called a \textit{path} of length $L$,
and the set of all paths of length $L$ is denoted by $\D_L$. 
One can illustrate a path graphically by drawing all points
$(i,\sigma_i)$ and connecting adjacent points by straight lines.
An example for a path with $L=9$ is given in Figure \ref{fig path}.
The condition $\sigma_i\sigma_{i+1}=0$ requires that the paths 
consist of a certain number of non-overlapping triangles
(this condition comes directly from the condition in the two-dimensional 
hard hexagon model requiring that no two particles can be on adjacent
sites).
To each path $\sigma=(\sigma_0,\ldots,\sigma_L)$ one may assign an energy 
$E(\sigma)$ by summing up the positions of the peaks, that is
\begin{equation*}
E(\sigma)=\sum_{j=1}^L j \sigma_j.
\end{equation*}
The energy of the path in Figure \ref{fig path}
is $E(\sigma)=1+5+8=14$. The generating function of paths of length $L$ 
which is also called a one dimensional configuration sum is defined as 
\begin{equation}\label{path}
X(L)=\sum_{\sigma\in \D_L} q^{E(\sigma)}.
\end{equation}
The path picture immediately implies that $X(L)$
satisfies the following initial conditions and recurrence
which completely specify it
\begin{equation}
\label{rec}
\begin{split}
X(0)&=X(1)=1\\
X(L)&=X(L-1)+q^{L-1}X(L-2).
\end{split}
\end{equation}

The aim is to find explicit expressions for $X(L)$.
We will describe two ways to obtain such an expression
which will be related to the two sides of \eqref{bf1}.

\subsection{Fermionic formula}
\label{subsec:fermi}

\begin{figure}
\begin{picture}(200,50)(0,0)
\put(10,20){\line(1,1){20}}
\put(50,20){\line(-1,1){20}}
\put(50,20){\line(1,1){20}}
\put(90,20){\line(-1,1){20}}
\put(90,20){\line(1,1){20}}
\put(130,20){\line(-1,1){20}}
\put(8,9){0}
\put(28,9){1}
\put(48,9){2}
\put(68,9){3}
\put(88,9){4}
\put(108,9){5}
\put(128,9){6}
\put(148,9){7}
\put(168,9){8}
\put(188,9){9}
\put(10,20){\line(1,0){190}}
\put(10,20){\circle*{4}}
\put(30,40){\circle*{4}}
\put(50,20){\circle*{4}}
\put(70,40){\circle*{4}}
\put(90,20){\circle*{4}}
\put(110,40){\circle*{4}}
\put(130,20){\circle*{4}}
\put(150,20){\circle*{4}}
\put(170,20){\circle*{4}}
\put(190,20){\circle*{4}}
\end{picture}
\caption{Ground state path with three particles}
\label{fig:ground}
\end{figure}
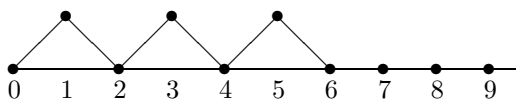

Interpret each peak in a path as a particle, that is, there
is a particle at site $i$ if $\sigma_i=1$.
Fix the number of particles to be $n$. The ground state path $\sg$ with
minimal energy is the path with particles at positions
$1,3,\ldots,2n-1$. The energy of $\sg$ is
$E(\sg)=1+3+\cdots+2n-1=n^2$. An example of the ground state path with $3$
particles is shown in Figure \ref{fig:ground}.
All other paths with $n$ particles can be obtained from $\sg$
by moving particles to the right in such a way that the particles
never overtake each other. If the length of the path is $L$,
the rightmost particle can move at most $L-2n$ positions to the right.
Hence the paths of length $L$ with $n$ particles are in one-to-one 
correspondence with partitions with at most $n$ parts not exceeding
$L-2n$. If a path $\sigma$ corresponds to partition $\la$ then
its energy is $E(\sigma)=E(\sg)+|\la|=n^2+|\la|$. The generating function
of partitions with at most $n$ parts not exceeding
$m$ is the $q$-binomial coefficient defined as
\begin{equation*}
\qbin{m+n}{n}=\frac{(q)_{m+n}}{(q)_m(q)_n}
\end{equation*}
for $m,n\in\N$ and zero otherwise. This implies the following
explicit expression for $X(L)$
\begin{equation}\label{RR fermi}
X(L)=\sum_{n=0}^\infty q^{n^2} \qbin{L-n}{n}.
\end{equation}
Because of its interpretation in terms of non-overlapping
particles this expression is called a quasiparticle or fermionic formula.

\subsection{Bosonic formula}
\label{subsec:bose}

As opposed to \eqref{RR fermi} there is another expression
for $X(L)$ given by
\begin{equation}\label{RR bose}
X(L)=\sum_{j=-\infty}^\infty (-1)^j q^{\frac{j}{2}(5j+1)}
 \qbin{L}{\lfloor \frac{L-5j}{2}\rfloor}
\end{equation}
where $\lfloor x\rfloor$ is the largest integer not exceeding $x$.
It can be proven by showing that the right-hand side
satisfies the recurrence and initial condition \eqref{rec}.

\begin{figure}
\begin{picture}(210,80)(-10,0)
\put(10,60){\line(1,1){20}}
\put(30,80){\line(1,-1){20}}
\put(50,60){\line(1,-1){20}}
\put(70,40){\line(1,1){20}}
\put(90,60){\line(1,1){20}}
\put(110,80){\line(1,-1){20}}
\put(130,60){\line(1,-1){20}}
\put(150,40){\line(1,-1){20}}
\put(170,20){\line(1,1){20}}
\put(8,9){0}
\put(28,9){1}
\put(48,9){2}
\put(68,9){3}
\put(88,9){4}
\put(108,9){5}
\put(128,9){6}
\put(148,9){7}
\put(168,9){8}
\put(188,9){9}
\put(-3,17){1}
\put(-3,37){2}
\put(-3,57){3}
\put(-3,77){4}
\put(10,20){\line(1,0){180}}
\put(10,80){\line(1,0){180}}
\put(10,20){\line(0,1){60}}
\put(190,20){\line(0,1){60}}
\put(10,60){\circle*{4}}
\put(30,80){\circle*{4}}
\put(50,60){\circle*{4}}
\put(70,40){\circle*{4}}
\put(90,60){\circle*{4}}
\put(110,80){\circle*{4}}
\put(130,60){\circle*{4}}
\put(150,40){\circle*{4}}
\put(170,20){\circle*{4}}
\put(190,40){\circle*{4}}
\end{picture}
\caption{Path associated to the path in Figure \ref{fig path}
under the correspondence \eqref{equiv}}
\label{fig:newpath}
\end{figure}
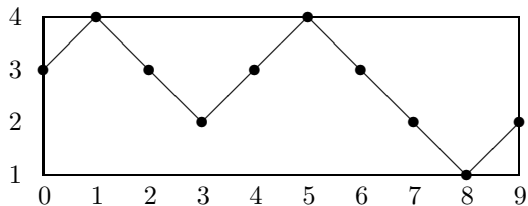

However, formula \eqref{RR bose} can also be interpreted
in terms of paths. For this purpose the above path picture needs
to be slightly reformulated. The state diagram of the paths
in the hard hexagon model is
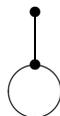
\begin{equation*}
\begin{picture}(20,40)(-10,0)
\put(0,20){\line(0,1){20}}
\put(0,10){\circle{20}}
\put(0,20){\circle*{4}}
\put(0,40){\circle*{4}}
\end{picture}
\end{equation*}
This diagram is to be understood as follows. If the system is in the
bottom state at position $i$, that is $\sigma_i=0$, then at position 
$i+1$ it can either be in the bottom state again so that $\sigma_{i+1}=0$
as indicated by the loop or it can be in the top state which means 
$\sigma_{i+1}=1$ as indicated by the line.
On the other hand, if the system is in the top state at position $i$
then at the next position it has to be in the bottom state, meaning
that $\sigma_i=1$ implies $\sigma_{i+1}=0$. This corresponds to
condition on paths that $\sigma_i\sigma_{i+1}=0$ for all $i$.
Up to a $\Z_2$ symmetry this state diagram is isomorphic to
\begin{equation}
\label{equiv}
\raisebox{-.3cm}{
\begin{picture}(20,40)(-10,0)
\put(0,20){\line(0,1){20}}
\put(0,10){\circle{20}}
\put(0,20){\circle*{4}}
\put(0,40){\circle*{4}}
\end{picture}
}
\quad\cong\quad
\raisebox{-1cm}{
\begin{picture}(4,60)(-2,0)
\put(0,0){\line(0,1){60}}
\put(0,0){\circle*{4}}
\put(0,20){\circle*{4}}
\put(0,40){\circle*{4}}
\put(0,60){\circle*{4}}
\end{picture}}
\;\;\Big{/}\;\; \Z_2.
\end{equation}
Under this correspondence the paths of the hard hexagon model
become paths in a strip of height 4 with steps going up or down
at each position. For example, with the convention that the paths start at
height three, the path in Figure \ref{fig path} becomes
the path in Figure \ref{fig:newpath}.

Let us now consider the following set of paths. Let
$p=(p_1,p_2,\ldots,p_L)$ be a sequence of $1$'s and $2$'s
and let $\la=(\la_1,\la_2)$ be the content of $p$ where $\la_i$ is
the number of $i$'s in $p$. Denote by $\P_{L,\la}$ the set of all 
$p=(p_1,\ldots,p_L)$ with content $\la$. Another way to represent
$p$ is by height variables $\sigma=(\sigma_0,\sigma_1,\ldots,\sigma_L)$
where by convention $\sigma_0=3$ and $\sigma_i=\sigma_{i-1}+1$ if
$p_i=1$ and $\sigma_i=\sigma_{i-1}-1$ if $p_i=2$ for $1\le i\le L$.
We will use $\sigma$ and $p$ interchangeably.
We are actually mostly concerned here with the set $\P_L:=\P_{L,\la}$
where $\la=(\lfloor \frac{L}{2}\rfloor,\lfloor \frac{L+1}{2}\rfloor)$. 

Consider the following subsets of $\P_L$.
Let $\Pbar_L$ be the subset of $\P_L$ consisting of
all paths in the strip of height four, that is all $\sigma\in\P_L$
such that $1\le \sigma_i\le 4$ for all $1\le i\le L$.
Furthermore, let $\P_L^{\uparrow,j}$ be the set of all paths 
$\sigma\in\P_L$ such that there exist indices $1\le i_1<i_2<\cdots<i_j\le L$
so that $\sigma_{i_1},\sigma_{i_3},\sigma_{i_5},\ldots$ are greater
than four and $\sigma_{i_2},\sigma_{i_4},\sigma_{i_6},\ldots$ are less than
one. Similarly let $\P_L^{\downarrow,j}$ be the set of all paths 
$\sigma\in\P_L$ such that there exist indices 
$1\le i_1<i_2<\cdots<i_j\le L$ such that 
$\sigma_{i_1},\sigma_{i_3},\sigma_{i_5},\ldots$ are less than one
and $\sigma_{i_2},\sigma_{i_4},\sigma_{i_6},\ldots$ are greater than
four. By inclusion-exclusion we have
\begin{equation*}
\Pbar_L=\Bigl(\P_L \cup \bigcup_{j\ge 1} (
\P_L^{\uparrow,2j} \cup \P_L^{\downarrow,2j}) \Bigr) \setminus
\bigcup_{j\ge 1} (\P_L^{\uparrow,2j-1} \cup \P_L^{\downarrow,2j-1}).
\end{equation*}

Now define the energy function $E$ for $\sigma\in\P_L$ as
\begin{equation*}
E(\sigma)=\sum_{i=1}^{L-1} i\cdot h(\sigma_{i-1},\sigma_i,\sigma_{i+1})
\end{equation*}
where $h$ is the local energy function given by
\begin{equation*}
h(\sigma_{i-1},\sigma_i,\sigma_{i+1})=
\begin{cases} 1 & \text{if $\sigma_{i-1}=\sigma_i-1=\sigma_{i+1}$ and
 $\sigma_i>3$,}\\
 & \text{or $\sigma_{i-1}=\sigma_i+1=\sigma_{i+1}$ and
 $\sigma_i<2$,}\\
0&\text{else.}
\end{cases}
\end{equation*}
The term $q^{\frac{j}{2}(5j+1)}\qbin{L}{\lfloor \frac{L-5j}{2} \rfloor}$
in \eqref{RR bose} can then be interpreted as the generating function
of $\P_L^{\downarrow,j}$ given by $\sum_{\sigma\in\P_L^{\downarrow,j}}
q^{E(\sigma)}$ if $j>0$ and the generating function of $\P_L^{\uparrow,j}$ 
given by $\sum_{\sigma\in\P_L^{\uparrow,j}} q^{E(\sigma)}$ if $j<0$. 
The term $j=0$ in \eqref{RR bose} is the generating function
of $\P_L$.

In section \ref{sec:bose} we will encounter more general 
inclusion-exclusion arguments in terms of operations on crystals.

\subsection{Identities}

Equations \eqref{RR fermi} and \eqref{RR bose} yield the
following polynomial identity
\begin{equation*}
\sum_{n=0}^\infty q^{n^2} \qbin{L-n}{n}=
\sum_{j=-\infty}^\infty (-1)^j q^{\frac{j}{2}(5j+1)}
 \qbin{L}{\lfloor \frac{L-5j}{2}\rfloor}.
\end{equation*}
Using $\lim_{L\to\infty} \qbins{L}{m}=\frac{1}{(q)_m}$ this identity
implies \eqref{bf1} in the limit $L\to\infty$.

All discussion so far focused on \eqref{bf1}. The second
Rogers--Ramanujan identity \eqref{bf2} can in fact be treated
in a very similar fashion. Define the set $\D'_L$ analogous to
$\D_L$ with the only difference that now $\sigma_0=1$ instead of 0.
The generating function is defined to be 
$X'(L)=\sum_{\sigma\in \D'_L} q^{E(\sigma)}$.
Analogous arguments to those in sections \ref{subsec:fermi} and 
\ref{subsec:bose} yield the following two expressions for $X'(L)$
\begin{equation*}
X'(L)=\sum_{n=0}^\infty q^{n(n+1)} \qbin{L-n-1}{n}
 =\sum_{j=-\infty}^\infty (-1)^j q^{\frac{j}{2}(5j+3)}
 \qbin{L}{\lfloor \frac{L-5j-1}{2}\rfloor}.
\end{equation*}
In the limit $L\to\infty$ this polynomial identity implies \eqref{bf2}.

\section{Crystal bases}
\label{sec:crystal}

In this section we review the main results about crystal bases
and explain how they provide a general setting for the definition
of one-dimensional configuration sums.

The quantized universal enveloping algebra $U_q(\gggg)$
associated with a symmetrizable Kac--Moody Lie algebra $\gggg$
was introduced independently by Drinfeld \cite{D} and Jimbo \cite{J}
in their study of two dimensional solvable lattice models in
statistical mechanics. The parameter $q$ corresponds 
to the temperature of the underlying model.
Kashiwara \cite{K} showed that at zero temperature or $q=0$ the
representations of $U_q(\gggg)$ have bases,
which he coined crystal bases, with a beautiful combinatorial structure
and favorable properties such as uniqueness and stability under tensor
products. The existence and uniqueness of crystal bases for integrable
highest weight modules for an arbitrary symmetrizable Kac--Moody algebra 
was given in \cite{K1}.

\subsection{Axiomatic definition of crystals}

Let $\gggg$ be a symmetrizable Kac-Moody algebra, $P$ the weight
lattice, $I$ the index set for the vertices of the Dynkin diagram
of $\gggg$, $\{\alpha_i\in P \mid i\in I \}$ the simple roots, and
$\{h_i\in P^*=\mathrm{Hom}_{\Z}(P,\Z)\mid i\in I \}$ the simple coroots.
Let $U_q(\gggg)$ be the quantized universal enveloping algebra of
$\gggg$. A $U_q(\gggg)$-crystal is a nonempty set $B$
equipped with maps $\wt:B\rightarrow P$ and
$e_i,f_i:B\rightarrow B\cup\{\emptyset\}$ for all $i\in I$,
satisfying
\begin{align}
&f_i(b)=b' \Leftrightarrow e_i(b')=b
\text{ if $b,b'\in B$} \\
&\wt(f_i(b))=\wt(b)-\alpha_i \text{ if $f_i(b)\in B$} \\
&\inner{h_i}{\wt(b)}=\varphi_i(b)-\epsilon_i(b)
\end{align}
where $\langle \cdot,\cdot\rangle$ is the natural pairing.
Here for $b\in B$,
\begin{equation} \label{eq:eps phi def}
\begin{split}
\epsilon_i(b)&= \max\{n\ge0\mid e_i^n(b)\not=\emptyset \} \\
\varphi_i(b) &= \max\{n\ge0\mid f_i^n(b)\not=\emptyset \}.
\end{split}
\end{equation}
(It is assumed that $\varphi_i(b),\epsilon_i(b)<\infty$ for all
$i\in I$ and $b\in B$.) 

A $U_q(\gggg)$-crystal $B$ can be viewed as a directed edge-colored graph 
(the crystal graph) whose vertices are the elements of $B$, with a directed 
edge from $b$ to $b'$ labeled $i\in I$, if and only if $f_i(b)=b'$.
The element $b$ is said to be highest weight if $e_i(b)=\emptyset$ 
for all $i\in I$.
Let $K$ be a subset of $I$. Then the $K$-component of the 
crystal $B$ is the graph obtained by only considering
edges colored by $i\in K$.

We also define the crystal reflection operator $s_i:B\rightarrow B$ by
\begin{equation*}
s_i(b) = \begin{cases}
  f_i^{\varphi_i(b)-\epsilon_i(b)}(b) & 
   \text{if $\varphi_i(b)>\epsilon_i(b)$} \\
  b & \text{if $\varphi_i(b)=\epsilon_i(b)$} \\
  e_i^{\epsilon_i(b)-\varphi_i(b)}(b) & 
   \text{if $\varphi_i(b)<\epsilon_i(b)$}.
	\end{cases}
\end{equation*}
It is obvious that $s_i$ is an involution.

\subsection{Tensor products of crystals}

Given two crystals $B$ and $B'$, there is also a crystal
obtained by taking the tensor product $B\otimes B'$.
As a set $B\otimes B'$ is just given by the Cartesian product
of the sets $B$ and $B'$. The weight function $\wt$ for
$b\otimes b'\in B\otimes B'$ is $\wt(b\otimes b')=\wt(b)+\wt(b')$
and the raising and lowering operators $e_i$ and $f_i$ act as follows
\begin{equation}
\label{tensor}
\begin{split}
e_i(b\otimes b')&=\begin{cases}
 e_i b \otimes b' & \text{if $\epsilon_i(b)>\phi_i(b')$,}\\
 b \otimes e_i b' & \text{otherwise,}
\end{cases}\\
f_i(b\otimes b')&=\begin{cases}
 f_i b \otimes b' & \text{if $\epsilon_i(b)\ge\phi_i(b')$,}\\
 b \otimes f_i b' & \text{otherwise.}
\end{cases}
\end{split}
\end{equation}
The reader is warned that this convention is different from
Kashiwara's convention. The order of the tensor factors is
interchanged.

\subsection{Finite and infinite crystals}

Let us fix some notation. From now on
let $\gggg$ be a simple complex Lie algebra and $\gggh$ be the
associated untwisted affine algebra. That is, let
$\gggh'=\gggg\otimes\CC[t,t^{-1}] \oplus \CC c$ be the central
extension of the loop algebra of $\gggg$ and
$\gggh=\gggh'\oplus\CC d$ where $d$ is the degree derivation. Let
$J$ (resp. $I=J\cup\{0\}$) be a set indexing the vertices of the
Dynkin diagram of $\gggg$ (resp. $\gggh$ and $\gggh'$).
A classical weight is a weight with respect to the algebra $\gggg$.
Denote by $U_q(\gggg)$, $U'_q(\gggh)$, and $U_q(\gggh)$ the
quantized universal enveloping algebras of $\gggg$, $\gggh'$, and
$\gggh$ respectively.

There are two main categories of crystals. The first one
contains the crystal bases of irreducible integrable
$U_q(\gggh)$-modules $V(\La)$ of highest weight $\La\in P^+$
where $P^+$ denotes the set of dominant weights of $\gggh$.
These crystals are infinite-dimensional. The second one
contains finite-dimensional crystals corresponding to 
$U_q'(\gggh)$. Since the set $B$
is finite in this case these crystals are called finite 
crystals. The level of a finite crystal $B$ is defined as
\begin{equation*}
\lev B = \min\{ \langle c,\epsilon(b) \rangle\mid b\in B\}
\end{equation*}
where $\epsilon(b)=\sum_{i\in I} \epsilon_i(b) \La_i$ and
$\{\La_i\mid i\in I\}$ is the $\Z$-basis of the weight lattice of
$\gggh'$.
A $U'_q(\gggh)$-crystal can be viewed as a $U_q(\gggg)$-crystal
by restricting to the $J$-component.

The crystal of each integrable highest weight $U_q(\gggh)$-module
can be realized in terms of a semi-infinite tensor product of perfect 
crystals \cite{KKMMNN,KKMMNN2,KKM}.
Perfect crystals are finite crystals with some additional
properties (for details see \cite[Definition 4.6.1]{KKMMNN}).
At least one perfect crystal for each integrable
$U_q(\gggh)$-module for $\gggh=A_n^{(1)}, B_n^{(1)},
C_n^{(1)}, D_n^{(1)}, A_n^{(2)}$ and $D_{n+1}^{(2)}$ was given
in \cite{KKMMNN2}.

Kashiwara and Nakashima \cite{KN} constructed the finite
$U_q(\gggg)$-crystals for $\gggg=A_n$, $B_n$, $C_n$ and $D_n$ 
explicitly in terms of tableaux. The cases $\gggg=A_n$ and $C_n$ 
are discussed in more detail in the next subsection.

\subsection{Finite crystals of type $A_n$ and $C_n$}
\label{sec:finite AC}

The finite crystals associated with $\gggg=A_n$ and $C_n$ are presented
in more detail since they will be our main examples throughout
the paper. Let $\Lambda_i$ for $i\in J$ be the fundamental weights of $\gggg$.
For later purposes it will be convenient to give the root
and weight structure of $A_n$ and $C_n$ explicitly.
Let $(\cdot|\cdot)$ be the standard bilinear form normalized such that
$(\alpha_i\mid\alpha_i)=2$ for the long roots $\alpha_i$.
\begin{example}
\label{ex:roots}
For either $\gggg=A_n$ or $\gggg=C_n$ the index set for the simple roots
is $J=\{1,2,\dotsc,n\}$. 

For $\gggg=A_n$ the weight lattice is embedded in the subspace of
$\R^{n+1}$ orthogonal to the vector $e=\sum_{j=1}^{n+1} \epsilon_j$,
where $\epsilon_i$ is the $i$-th standard basis vector of $\R^{n+1}$.
The simple roots are $\alpha_i=\epsilon_i-\epsilon_{i+1}$ for all
$i\in J$. The fundamental weights are $\La_i=\sum_{j=1}^i \epsilon_j
-\frac{i}{n+1}e$ for $i\in J$. The half-sum of positive roots
is $\rho=(n,n-1,\dotsc,1,0)-\frac{n}{2}e$.
The scalar product is the standard one:
$(\alpha|\beta)=\alpha\cdot\beta$.

For $\gggg=C_n$ the simple roots
are given by the short roots $\alpha_i=\epsilon_i-\epsilon_{i+1}$
for $1\le i<n$ and the long root $\alpha_n=2\epsilon_n$ where 
$\epsilon_i$ is the $i$-th unit vector in $\R^n$. The fundamental
weights are $\La_i=\epsilon_1+\cdots+\epsilon_i$ for $i\in J$.
Half the sum of all positive roots is $\rho=(n,n-1,\ldots,1)$.
The bilinear form is given by $(\alpha|\beta)=\alpha\cdot \beta/2$.

For both type $A_n$ and $C_n$ we will identify dominant weights,
that is weights of the form $\La=\La_{k_1}+\cdots+\La_{k_m}$,
with partitions. 
A partition $\la=(\la_1,\ldots,\la_{n+1})$ (with $\la_{n+1}=0$ for
type $C_n$) corresponds to the weight 
$\La=\sum_{i\in J}(\la_i-\la_{i+1})\La_i$.
\end{example}

For each dominant weight $\Lambda$ there is a finite classical crystal
$B(\Lambda)$ \cite{KN}. The crystals $B(\Lambda_1)$ for type $A_n$ and 
$C_n$ are given in Table \ref{table: classical crystals} where the arrow
$\stackrel{i}{\longrightarrow}$ stands for $f_i$.
\begin{table}
\begin{tabular}{|c|l|}
\hline
%
$A_n$ &
\raisebox{-0.3cm}{\scalebox{0.7}{
\begin{picture}(250,38)(-10,-19)
\BText(0,0){1}
\LongArrow(10,0)(40,0)
\BText(50,0){2}
\LongArrow(60,0)(90,0)
\BText(100,0){3}
\LongArrow(110,0)(140,0)
\Text(160,0)[]{$\cdots$}
\LongArrow(175,0)(205,0)
\BText(220,0){n+1}
\PText(25,2)(0)[b]{1}
\PText(75,2)(0)[b]{2}
\PText(125,2)(0)[b]{3}
\PText(190,2)(0)[b]{n}
\end{picture}
}}
\\ \hline
%
$C_n$ &
\raisebox{-0.3cm}{\scalebox{0.7}{
\begin{picture}(350,38)(-10,-19)
\BText(0,0){1}
\LongArrow(10,0)(30,0)
\BText(40,0){2}
\LongArrow(50,0)(70,0)
\Text(85,0)[]{$\cdots$}
\LongArrow(95,0)(115,0)
\BText(125,0){n}
\LongArrow(135,0)(155,0)
\BBoxc(165,0)(13,13)
\Text(165,0)[]{\footnotesize$\overline{\mbox{n}}$}
\LongArrow(175,0)(195,0)
\Text(210,0)[]{$\cdots$}
\LongArrow(220,0)(240,0)
\BBoxc(250,0)(13,13)
\Text(250,0)[]{\footnotesize$\overline{\mbox{2}}$}
\LongArrow(260,0)(280,0)
\BBoxc(290,0)(13,13)
\Text(290,0)[]{\footnotesize$\overline{\mbox{1}}$}
\PText(20,2)(0)[b]{1}
\PText(60,2)(0)[b]{2}
\PText(105,2)(0)[b]{n-1}
\PText(145,2)(0)[b]{n}
\PText(185,2)(0)[b]{n-1}
\PText(230,2)(0)[b]{2}
\PText(270,2)(0)[b]{1}
\end{picture}
}}
\\ \hline
\end{tabular}\vspace{4mm}
\caption{\label{table: classical crystals}Classical crystals $B(\La_1)$}
\end{table}
The finite crystals $B(\La_k)$ for $k\in J$ can be obtained in
the following way. Let $u_{\La_k}=
\raisebox{-0.14cm}{\scalebox{0.8}{\begin{picture}(14,14)(0,-7)
\BBoxc(7,0)(14,14)\Text(7,0)[]{$k$}\end{picture}}}
\otimes \cdots \otimes
\raisebox{-0.14cm}{\scalebox{0.8}{\begin{picture}(14,14)(0,-7)
\BBoxc(7,0)(14,14)\Text(7,0)[]{$2$}\end{picture}}}
\otimes
\raisebox{-0.14cm}{\scalebox{0.8}{\begin{picture}(14,14)(0,-7)
\BBoxc(7,0)(14,14)\Text(7,0)[]{$1$}\end{picture}}}$
be the unique highest weight vector of weight $\La_k$
in $B(\La_1)^{\otimes k}$. Then $B(\La_k)$ is the connected
component of $B(\La_1)^{\otimes k}$ containing $u_{\La_k}$.
In \cite{KN} the elements in this connected component were
identified with certain one-column tableaux.
The finite crystal $B(\La)$ for a dominant weight $\La=\La_{k_1}+\cdots+
\La_{k_m}$ with $k_1\ge k_2\ge \cdots \ge k_m$ is isomorphic to the
connected component in $B(\La_{k_1})\otimes \cdots \otimes B(\La_{k_m})$
which contains the highest weight element $u_{\La_{k_1}}\otimes \cdots
\otimes u_{\La_{k_m}}$.
Combining the two embeddings it follows that $B(\La)$ is isomorphic
to a certain connected component in $B(\La_1)^{\otimes M}$ where
$M=k_1+\cdots+k_m$.

For type $A_n$ the elements in $B(\La)$ can be identified
with semi-standard Young tableaux of shape $\la$ where $\la$
is the partition corresponding to the weight $\La$.
The paths encountered in section \ref{sec:HHM} are sequences of
1's and 2's. Viewed as single box Young tableaux over the
alphabet $\{1,2\}$, these are exactly the elements of the 
crystal $B(\La_1)$ of type $A_1$.

\subsection{Simple crystals}

Simple crystals were introduced by Akasaka and Kashiwara \cite{AK}.
As will be explained in the next section they have an
isomorphism and energy function which are required for the
definition of the one-dimensional configuration sums.

Let $\Wh$ be the Weyl group of the affine Kac--Moody algebra $\gggh$
generated by the simple reflections $r_i$ for $i\in I$ defined as
\begin{equation}
\label{reflections}
r_i(\beta)=\beta - \langle h_i,\beta\rangle \alpha_i
\end{equation}
where $\beta$ is a root.
Let $B$ be the crystal graph of an integrable $U_q(\gggh)$-module.
Say that $b\in B$ is an extremal vector of weight $\Lambda\in P$
provided that $\wt(b)=\Lambda$ and there exists a family of elements
$\{b_w\mid w\in \Wh \}\subset B$ such that
\begin{enumerate}
\item $b_w=b$ for $w=e$.
\item If $\inner{h_i}{w\La}\ge 0$ then $e_i(b)=0$ and
$f_i^{\inner{h_i}{w\La}}(b_w)=b_{r_i w}$.
\item If $\inner{h_i}{w\La}\le 0$ then $f_i(b)=0$ and
$e_i^{\inner{h_i}{w\La}}(b_w)=b_{r_i w}$.
\end{enumerate}

\begin{definition}
Say that a $U'_q(\gggh)$-crystal $B$ is \textit{simple} if
\begin{enumerate}
\item $B$ is the crystal base of a finite-dimensional integrable
$U'_q(\gggh)$-module.
\item There is a dominant weight $\Lambda$ with respect to the weight
lattice of the classical algebra $\gggg$ such that $B$ has a
unique vector (denoted $u(B)$) of weight $\Lambda$, and the weight of
any extremal vector of $B$ is contained in $W\Lambda$.
Here $W$ is the Weyl group corresponding to the classical algebra $\gggg$.
\end{enumerate}
\end{definition}

\begin{theorem} \label{thm:simple} \cite{AK}
\begin{enumerate}
\item Simple crystals are connected as graphs.
\item The tensor product of simple crystals is simple.
\end{enumerate}
\end{theorem}

Let $B$ be a simple $U'_q(\gggh)$-crystal, equipped with a
function $D=D_B:B\rightarrow\Z$, called its intrinsic energy, which
is required to be constant on $J$-components and defined up to a
global additive constant. Call the pair $(B,D)$ a graded simple
$U'_q(\gggh)$-crystal. We normalize the intrinsic energy function
by the requirement that
\begin{equation*}
D_B(u(B))=0.
\end{equation*}

\subsection{Finite dimensional affine crystals}

Recently, new families of crystals of the finite
dimensional representations of $U_q(\gggh)$ were conjectured 
\cite{HKOTY,K3} where $\gggh$ is a untwisted affine Lie algebra.

\begin{conjecture}\cite{HKOTY,K3}\label{conj:Brs}
For each $r\in J$ and $s\ge1$, there exists an irreducible
finite-dimensional integrable $U'_q(\gggh)$-module $W^{(r)}_s$
with simple crystal base $B^{r,s}$ generated a unique extremal vector
$u(B^{r,s})$ of weight $s\La_r$, and a prescribed
$U_q(\gggg)$-crystal decomposition of the form $B^{r,s} \cong
B(s\La_r)\oplus B$, where $B$ is a direct sum of
$U_q(\gggg)$-crystals of the form $B(\La)$ where $\La$ is a classical
dominant weight and $s \La_r \vartriangleright \la$. 
Here $\La'\trianglerighteq\La$ if and only if
$\La'-\La\in\bigoplus_{i\in J} \N \alpha_i$.
Moreover there is a prescribed intrinsic energy function
$D=D_{B^{r,s}}:B^{r,s}\rightarrow \Z$, that is constant on
$J$-components, such that $0=D(u)>D(b)$ where $u$ is the
$J$-highest weight vector of weight $s\La_r$ in $B^{r,s}$, and
$b$ is any element not in the $J$-component of $u$.
\end{conjecture}

\subsection{Combinatorial $R$-matrix and energy function}

Suppose $B_1$ and $B_2$ are simple $U'_q(\gggg)$-crystals. Then
there is a unique isomorphism of $U'_q(\gggg)$-crystals
$\sigma:B_2\otimes B_1\rightarrow B_1\otimes B_2$.
There is also a function $H=H_{B_2,B_1}:B_2\otimes B_1\rightarrow\Z$,
called local energy function which is unique 
up to global additive constant, such that for all 
$b_2\otimes b_1\in B_2\otimes B_1$ with $b_1'\otimes b_2'=
\sigma(b_2 \otimes b_1)$
\begin{equation} \label{eq:local energy}
  H(e_i(b_2\otimes b_1))=
  H(b_2\otimes b_1)+
  \begin{cases}
    -1 & \text{if $i=0$, $\epsilon_0(b_2)>\varphi_0(b_1)$,
    $\epsilon_0(b_1')>\varphi_0(b_2')$} \\
    1 & \text{if $i=0$, $\epsilon_0(b_2)\le\varphi_0(b_1)$,
    $\epsilon_0(b_1')\le \varphi_0(b_2')$} \\
    0 & \text{otherwise.}
  \end{cases}
\end{equation}
The pair $(\sigma,H)$ is called the combinatorial $R$-matrix.
It is convenient to normalize the local energy function $H$ by
requiring that
\begin{equation} \label{eq:extremal H}
  H(u(B_2)\otimes u(B_1)) = 0.
\end{equation}
With this convention it follows by definition that
\begin{equation} \label{eq:HR=H}
H_{B_1,B_2} \circ \sigma_{B_2,B_1} = H_{B_2,B_1}
\end{equation}
as operators on $B_2\otimes B_1$.

Let $(B_j,D_j)$ be graded simple $U'_q(\gggh)$-crystals for $1\le
j\le L$ and set $B=B_L\otimes\dotsm\otimes B_1$.
Following \cite{NY} define the energy function
$E_B:B\rightarrow\Z$ by
\begin{equation}\label{eq:NY}
E_B = \sum_{1\le i<j\le L} H_i\sigma_{i+1}\sigma_{i+2}\dotsm\sigma_{j-1}
\end{equation}
where $H_i$ (resp. $\sigma_i$) is the local energy function (resp. 
isomorphism) acting on the $i$-th and $(i+1)$-st tensor factor.
By the normalization assumption \eqref{eq:extremal H} it follows
that
\begin{equation}\label{eq:extremal NY}
E_B(u(B))=0.
\end{equation}
As shown in \cite{OSS}, the intrinsic energy $D_B$ for the $L$-fold 
tensor product $B=B_L\otimes\dotsm\otimes B_1$ is given by
\begin{equation}\label{eq:energy}
  D_B = E_B + \sum_{j=1}^L D_j \sigma_1\sigma_2\dotsm\sigma_{j-1}
\end{equation}
where $D_j$ acts on the rightmost tensor factor which is
$B_j$.

The coenergy and intrinsic coenergy are defined as
\begin{equation*}
 \Ebar_B=-E_B \qquad\text{and}\qquad \Dbar_B=-D_B.
\end{equation*}

\subsection{One-dimensional configuration sums}

There are three different sets of paths that we consider.
Let $B=B_L\otimes\cdots\otimes B_1$ where all $B_i$ are simple crystals.
For a classical weight $\Lambda$ the set of unrestricted paths
is defined as
\begin{equation}
\label{p}
\P(B,\La)=\{ b\in B \mid \wt(b)=\La \}.
\end{equation}
For a dominant classical weight $\La$ the
set of classically restricted paths is
\begin{equation}
\label{p'}
\P'(B,\La)=\{ b\in B \mid \text{$\wt(b)=\La$ and $e_ib=0$ for all 
 $i\in J$} \}
\end{equation}
and the set of level-restricted paths for $\ell\in\N$ is
\begin{equation}
\label{pl}
\P^\ell(B,\La)=\{ b\in B \mid \text{$\wt(b)=\La$, $e_ib=0$ for all 
 $i\in J$ and $e_0^{\ell+1}b=0$}\}.
\end{equation}

The corresponding one-dimensional configuration sums are
the generating functions of these sets of paths with
energy/coenergy statistics.
The one-dimensional sums
\begin{equation}
\label{super}
\begin{split}
S(B,\La)&=\sum_{b\in \P(B,\La)} q^{D_B(b)}\\
\Sbar(B,\La)&=\sum_{b\in \P(B,\La)} q^{\Dbar_B(b)}
\end{split}
\end{equation}
are called supernomials, whereas
\begin{equation}
\label{kostka}
\begin{split}
X(B,\La)&=\sum_{b\in \P'(B,\La)} q^{D_B(b)}\\
\Xbar(B,\La)&=\sum_{b\in \P'(B,\La)} q^{\Dbar_B(b)}
\end{split}
\end{equation}
are the classically restricted configuration sums or generalized Kostka
polynomials, and
\begin{equation}
\label{lkostka}
\begin{split}
X^\ell(B,\La)&=\sum_{b\in \P^\ell(B,\La)} q^{D_B(b)}\\
\Xbar^\ell(B,\La)&=\sum_{b\in \P^\ell(B,\La)} q^{\Dbar_B(b)}
\end{split}
\end{equation}
are the level-$\ell$ restricted configuration sums or
$\ell$-generalized Kostka polynomials.

The classically restricted configurations sums \eqref{kostka} are
graded tensor product multiplicities.
The level-restricted configuration sums \eqref{lkostka} are graded
level $\ell$ fusion coefficients. Let $B_{\La'}$ denote the crystal
corresponding to the affine irreducible highest weight representation
$V(\La')$. By the Verlinde formula \cite{V},
the fusion coefficient is the coefficient of $B_{\La'}$ of weight 
$\La' =\La+\ell\La_0$ in the decomposition of $B \otimes B_{\ell\La_0}$.
The affine highest weight vectors of weight $\La'$, whose number
is the above multiplicity, are the summands of $\Xbar^\ell(B,\La)$.

It will be shown in sections \ref{sec:bose} and \ref{sec:fermi}
that there are two different evaluations of $\Xbar(B,\La)$ and 
$\Xbar^\ell(B,\La)$ which give rise to $q$-identities.

\section{Bosonic evaluation}
\label{sec:bose}

Here we present the bosonic evaluation of $\Xbar(B,\Lambda)$ and 
$\Xbar^\ell(B,\Lambda)$ as defined in \eqref{kostka} and \eqref{lkostka}.
Similarly to the bosonic evaluation of $X(L)$ of section \ref{subsec:bose},
the bosonic evaluation of $\Xbar(B,\Lambda)$ and $\Xbar^\ell(B,\Lambda)$
can be obtained by inclusion-exclusion arguments as shown in \cite{SS}. 
We discuss the main ideas and techniques of sign-reversing involutions.

\subsection{Classically-restricted case}

Let $B=B_L\otimes \cdots \otimes B_1$ be a $U_q(\gggg)$-crystal 
and let $\La$ be a classical weight.
The Weyl group $W$ of $\gggg$ which is generated by the simple
reflections $r_i$ as in \eqref{reflections} with $i$ restricted to
$i\in J$.
The bosonic expression for the generating function of classically-restricted 
paths $\Xbar(B,\La)$ as defined in \eqref{kostka} is given by
\begin{equation}
\label{bose cr}
\Xbar(B,\La)=\sum_{\omega \in W} (-1)^\omega \; 
 \Sbar(B,\omega(\La+\rho)-\rho)
\end{equation}
where $(-1)^\omega$ is the sign of $\omega$ and $\rho$ is half the sum
of all positive roots.
This formula follows directly from Weyl's character formula.

As a warm-up for the level-restricted case, we would like to briefly
sketch for types $A_n$ and $C_n$ how \eqref{bose cr} can be derived via
a sign-reversing involution.
\begin{example}
\label{ex:Weyl}
For $\gggg=A_n$ the Weyl group $W$ is generated by the reflections
$r_1,\ldots,r_n$ which act on $\la\in\Z^{n+1}$ as follows
\begin{equation*}
r_i(\la)=(\la_1,\ldots,\la_{i+1},\la_i,\ldots,\la_{n+1}).
\end{equation*}
Hence on $\Z^{n+1}$ the Weyl group acts as the symmetric group $S_{n+1}$.

For $\gggg=C_n$ the Weyl group is generated by
\begin{equation*}
\begin{aligned}
r_i(\la)&=(\la_1,\ldots,\la_{i+1},\la_i,\ldots,\la_n) \qquad
 \text{for $1\le i<n$}\\
r_n(\la)&=(\la_1,\ldots,\la_{n-1},-\la_n)
\end{aligned}
\end{equation*}
so that $W$ on $\Z^n$ acts by all permutations and sign changes.
\end{example}

Set
\begin{equation*}
\mathcal{S}=\{(\omega,b)\in W\times B \mid \omega(\wt(b)+\rho)=\La+\rho\}
\end{equation*}
On the set $\mathcal{S}$ we define an involution
$\Phi:\mathcal{S}\to\mathcal{S}$ with the properties that
the fixed points are the pairs $(1,b)$ with $b\in \P'(B,\La)$ which
recall are the paths underlying $\Xbar(B,\La)$ (see \eqref{kostka}).
Furthermore, if $\Phi(\omega,b)=(\omega',b')$ with 
$(\omega,b)\neq (\omega',b')$ then $\omega$ and $\omega'$ have opposite signs.

Let $\mathcal{S}_i$ for $i\in J$ be the set of pairs 
$(\omega,b)\in \mathcal{S}$ such that $\epsilon_i(b)>0$. 
Define $\Phi_i:\mathcal{S}_i\to\mathcal{S}_i$ by
$\Phi_i(\omega,b)=(\omega r_i,s_i e_i(b))$.
Define the set $\mathcal{S}'=\mathcal{S}-\{(1,b)\mid b\in\P'(B,\La)\}$
so that $\mathcal{S}'=\bigcup_{i\in J} \mathcal{S}_i$. 
Then $\Phi$ is given by
\begin{equation*}
\Phi(\omega,b)=\begin{cases}
 (\omega,b) & \text{if $(\omega,b)\in \mathcal{S}\setminus \mathcal{S}'$}\\
 \Phi_i(\omega,b) & \text{if $(\omega,b)\in \mathcal{S}'$ and $i=v(\omega,b)$}
\end{cases}
\end{equation*}
where $v$ is some functions $v:\mathcal{S}'\to J$.
To show that $\Phi$ indeed exists and is an involution it 
needs to be shown that if $v(\omega,b)=i$ then
\begin{equation}
\label{v prop}
(\omega,b)\in\mathcal{S}_i \qquad \text{and} \qquad v(\Phi_i(\omega,b))=i.
\end{equation}
As mentioned in section \ref{sec:finite AC}, the crystal 
$B=B_L\otimes\cdots\otimes B_1$ of type $A_n$ or $C_n$ can be embedded into
$B(\La_1)^{\otimes M}$ for some $M$. Let $p=p_M\otimes\cdots\otimes p_1$ 
be the image of $b\in B$ under this embedding. Then $v$ can be defined
as follows. Let $k$ be minimal such that $\epsilon_i(p_k\otimes\cdots\otimes
p_1)>0$ for some $i\in J$. Then it is clear from 
Table \ref{table: classical crystals} that there is a unique $i$ satisfying 
$\epsilon_i(p_k\otimes\cdots\otimes p_1)>0$. Set $v(\omega,b)=i$.
It follows from equation \eqref{tensor} and 
Table \ref{table: classical crystals}
that the first $k$ tensor factors of $p$ stay invariant under 
$\Phi_i$ since there are no strings of length greater than one. 
This ensures \eqref{v prop}.
Hence inclusion-exclusion implies \eqref{bose cr}.

\subsection{Level-restricted case}

The bosonic expression for the level-restricted generating function
$\Xbar^\ell(B,\La)$ defined in \eqref{lkostka} can also
be found by a sign-reversing involution. The difference is that
one needs to consider elements $\omega$ in the affine Weyl group $\Wh$
which is generated by $r_i$ with $i\in I$ (rather than $i\in J$).

\begin{example}
For $\gggh=A_n^{(1)}$ the affine Weyl group $\Wh$ is generated
by the reflections $r_1,\ldots,r_n$ as in example \ref{ex:Weyl} and
\begin{equation*}
r_0(\la)=(\la_{n+1}+\ell+n+1,\la_2,\ldots,\la_n,\la_1-(\ell+n+1)).
\end{equation*}

For $\gggh=C_n^{(1)}$ the affine Weyl group $\Wh$ is generated by
the reflections $r_1,\ldots,r_n$ as in example \ref{ex:Weyl} and
\begin{equation*}
r_0(\la)=(-\la_1+2(\ell+n+1),\la_2,\ldots,\la_n).
\end{equation*}
\end{example}

The affine Weyl group is isomorphic to $\Wh \cong T \rtimes W$
where $T$ is the set of certain translations $t_\alpha$
indexed by $\alpha\in M$ for a particular set $M$
(for more details see \cite[Section 6]{Kac}).
Then it was shown in \cite{SS} that
\footnote{The arguments in \cite{SS} require that $B$ is a tensor product
of almost perfect crystals and that the energy function obeys certain
properties. For the examples of type $A_n^{(1)}$
with $B_i=B^{r_i,s_i}$ and $C_n^{(1)}$ with $B_i=B^{r_i,1}$, for which
we will consider fermionic formulas in the next section, these conditions
are all satisfied.}
\begin{multline}
\label{bose lr}
\Xbar^\ell(B,\La)=\sum_{\omega\in W} \sum_{\alpha\in M} (-1)^\omega 
 q^{\frac{1}{2}a_0 (\alpha|\alpha)(\ell+h^\vee)-a_0(\rho+\La|\alpha)}\\
 \qquad \times \Sbar(B,\omega(\La+\rho-(\ell+h^\vee)\alpha)-\rho)
\end{multline}
where $h^\vee$ is the dual Coxeter number of $\gggh$ and 
$a_0$ is the label of the zeroth node in the Dynkin diagram
of $\gggh$.

\begin{example}
Let us give \eqref{bose lr} more explicitly for
$\gggh=A_n^{(1)}$ and $C_n^{(1)}$.
In both cases the dual Coxeter number is $h^\vee=n+1$ and $a_0=1$.
For type $A_n^{(1)}$ the set $M$ is given by all $\beta\in \Z^{n+1}$ such that 
$|\beta|:=\beta_1+\cdots+\beta_{n+1}=0$ so that
\begin{multline*}
\Xbar^\ell(B,\la)=\sum_{\omega\in W}
 \sum_{\substack{\beta\in\Z^{n+1}\\ |\beta|=0}} (-1)^\omega
 q^{\frac{1}{2} \beta\cdot \beta(\ell+n+1)-(\rho+\la)\cdot \beta}\\
 \times \Sbar(B,\omega(\la+\rho-(\ell+n+1)\beta)-\rho)
\end{multline*}
where $W=S_{n+1}$ is the set of permutations.

For type $C_n^{(1)}$ we have $M=2\Z^n$. Hence
\begin{multline*}
\Xbar^\ell(B,\la)=\sum_{\omega\in W}
 \sum_{\beta\in 2\Z^n} (-1)^\omega
 q^{\frac{1}{4}\beta\cdot \beta (\ell+n+1)-\frac{1}{2}(\rho+\la)\cdot \beta}\\
 \times \Sbar(B,\omega(\la+\rho-(\ell+n+1)\beta)-\rho)
\end{multline*}
where the Weyl group $W$ in this case is the set of all permutations and 
sign changes.
The extra factor of $1/2$ in the exponent of $q$ comes from the
normalization of $(\cdot|\cdot)$ as alluded to in example \ref{ex:roots}.
\end{example}

The arguments in \cite{SS} involve a sign-reversing involution.
Similarly to the classically restricted case set
\begin{equation*}
\mathcal{S}=\{(\omega,b)\in \Wh\times B \mid \omega(\wt(b)+\rho)=\La+\rho\}.
\end{equation*}
For $i\in J$ define $\mathcal{S}_i$ and $\Phi_i$ as before. In addition,
let $\mathcal{S}_0$ be the subset of all pairs $(\omega,b)$ in 
$\mathcal{S}$ such that $\epsilon_0(b)>\ell$ and define
$\Phi_0(\omega,b)=(\omega r_0, s_0 e_0^{\ell+1} b)$.
One can find a sign-reversing involution $\Phi$ with fixed point set being
the set of level-$\ell$ restricted paths $\P^\ell(B,\La)$
by again specifying a function $v:\mathcal{S}-\{(1,b)\mid
b\in \P^\ell(B,\La)\}\to I$ satisfying \eqref{v prop}.
The existence of such a function $v$ was proven in \cite{SS}
for a large class of crystals.

\subsection{Supernomial coefficients}

The bosonic formulas \eqref{bose cr} and \eqref{bose lr}
involve the supernomial coefficients defined in \eqref{super}.
To obtain truely explicit expressions it is still necessary
to give formulas for the supernomial coefficients.
These are not yet known in general. Here we give a few
examples for which formulas exist.

\begin{example}
The supernomial coefficients for type $A_n^{(1)}$ for single columns were
given in \cite{HKKOTY,Ki}. Let $B=B^{\mu_L,1}\otimes \cdots \otimes 
B^{\mu_1,1}$ be the tensor product of crystals of type $A_n^{(1)}$
corresponding to the partition $\mu=(\mu_1,\ldots,\mu_L)$. Furthermore,
let $\la \in \N^{n+1}$ be a composition. Then
\begin{equation}
\label{super sc}
\Sbar(B,\la)=\sum_\nu \prod_{\substack{1\le a\le n\\ 1\le i\le \mu_1}}
 \qbin{\nu_i^{(a+1)}-\nu_{i+1}^{(a+1)}}{\nu_i^{(a)}-\nu_{i+1}^{(a+1)}}
\end{equation}
where the sum is over all sequences of partitions
$\nu=(\nu^{(1)},\ldots,\nu^{(n)})$ such that
\begin{equation*}
\begin{split}
&\emptyset=\nu^{(0)}\subset \nu^{(1)}\subset \cdots \subset \nu^{(n+1)}
=\mu^t\\
&\text{$\nu^{(a)}/\nu^{(a-1)}$ is a horizontal $\la_a$-strip}.
\end{split}
\end{equation*}
A horizontal $p$-strip is a skew shape with $p$ boxes such that
each column contains at most one box. 
For $\mu=(1^{|\la|})$ equation \eqref{super sc}
reduces to the $q$-multinomial coefficient
\begin{equation}
\label{multi}
\qbin{L}{\la_1,\ldots,\la_{n+1}}=\frac{(q)_L}{(q)_{\la_1}\cdots 
 (q)_{\la_{n+1}}} 
 \qquad \text{if $|\la|=L$}
\end{equation}
and zero otherwise.
\end{example}

\begin{example}
The supernomials for type $A_n^{(1)}$ for single rows were also
given in \cite{HKKOTY,Ki}.
Let $\mu=(\mu_1,\ldots,\mu_L)$ be a partition,
$B=B^{1,\mu_L}\otimes \cdots \otimes B^{1,\mu_1}$ and
$\la\in\N^{n+1}$ a composition. Then
\begin{equation}
\label{super sr}
\Sbar(B,\la)=\sum_\nu q^{\phi(\nu)} 
 \prod_{\substack{1\le a\le n\\ 1\le i\le \mu_1}}
 \qbin{\nu_i^{(a+1)}-\nu_{i+1}^{(a)}}{\nu_i^{(a)}-\nu_{i+1}^{(a)}}
\end{equation}
where the sum is over all sequences of partitions
$\nu=(\nu^{(1)},\ldots,\nu^{(n)})$ such that
\begin{equation*}
\begin{split}
&\emptyset=\nu^{(0)}\subset \nu^{(1)}\subset \cdots \subset \nu^{(n+1)}
=\mu^t\\
&\text{$|\nu^{(a)}|=\la_1+\cdots+\la_a$ for all $1\le a \le n$.}
\end{split}
\end{equation*}
Here $\phi(\nu)$ is defined as
\begin{equation*}
\phi(\nu)=\sum_{\substack{1\le a\le n\\ 1\le i\le \mu_1}}
  \nu_{i+1}^{(a)}(\nu_i^{(a+1)}-\nu_i^{(a)}).
\end{equation*}
As in the previous example this reduces to the $q$-multinomial coefficient
\eqref{multi} if $\mu=(1^{|\la|})$. For type $A_1$ this formula
coincides with that in \cite{SW}.
\end{example}

\begin{example}
The supernomials of type $C_n$ for single boxes can be obtained
by the following arguments. Let $B=(B^{1,1})^{\otimes L}$ and
$\la\in\Z^n$ with $\|\la\|:=|\la_1|+\cdots+|\la_n|\le L$. This means that
there are at least $\la_i$ letters $i$ if $\la_i\ge 0$ and
at least $\la_i$ letter $\overline{\imath}$ if $\la_i<0$. If $\|\la\|<L$ then
there have to be $(L-\|\la\|)/2$ pairs of barred and unbarred letters
in order to have weight $\la$. Hence we have
\begin{equation}
\Sbar(B,\la)=\sum_{\substack{\mu\in \N^n\\ 2|\mu|=L-\|\la\|}}
\qbin{L}{|\la_1|+\mu_1,\ldots,|\la_n|+\mu_n,\mu_1,\ldots,\mu_n}
\end{equation}
where $\qbin{L}{p_1,\ldots,p_k}$ is the $q$-multinomial as defined
in \eqref{multi}.
\end{example}

Explicit formulas also exist for level one cases for $B_n^{(1)}$, 
$D_n^{(1)}$ \cite{DJKMO} and $A_{2n-1}^{(2)}$, $A_{2n}^{(2)}$, 
$D_{n+1}^{(2)}$ \cite{KMOTU}.

\section{Fermionic evaluation}
\label{sec:fermi}

The derivation of fermionic evaluations of the classically-
and level-restricted configuration sums \eqref{kostka} and \eqref{lkostka}
is in general much more intricate than for the hard-hexagon model.
There exists a vast literature on conjectures and proofs of fermionic 
formulas, most of which deal with the case $\gggh=A_1^{(1)}$ or $A_n^{(1)}$.
A relatively complete list of references can be found in \cite{HKOTY}.
Fermionic formulas for all untwisted quantum affine algebras
were recently conjectured in \cite{HKOTY}. For type $A_n^{(1)}$ 
these are proven for the classically-restricted case in \cite{KSS}
and the level-restricted case in \cite{SSa}. For type $C_n^{(1)}$
the classically-restricted formulas are proven in \cite{OSS}. We will
present these results here and also derive the fermionic level-restricted
formulas for type $C_n^{(1)}$ in section \ref{sec:lev C}.

Interestingly, Kirillov and Reshetikhin \cite{KR} conjectured that
the coefficients of the decomposition of the representations of 
$U'_q(\gggh)$ naturally associated with multiples of the fundamental 
weights into direct sums of irreducible representations of $U_q(\gggg)$ 
are given by the fermionic formulas at $q=1$. 
Chari \cite{Chari} proved this conjecture for a single tensor factor
for $\gggg$ a simple Lie algebra of classical type and also for some
exceptional cases.

\subsection{Classically-restricted case}
\label{sec:classic}

We will state here the fermionic formulas conjectured in \cite{HKOTY}.
Let $\gggh=X_n^{(1)}$ with $X=A,B,C,D$ or $E$ for $n=6,7,8$ or $F$ for $n=4$
or $G$ for $n=2$.
Let $B=\bigotimes_{a=1}^n \bigotimes_{i\ge 1} (B^{a,i})^{\otimes L_i^{(a)}}$
where $L_i^{(a)}\in\N$ for all $i\ge 1$ and $1\le a\le n$ and only
finitely many $L_i^{(a)}$ are nonzero. Define the following polynomial
in $q$ depending on $B$ and a dominant weight $\La$
\begin{equation}\label{fermi kostka}
\Fbar(B,\La)=\sum_{\{m\}} q^{\cc(\{m\})} \prod_{a=1}^n \prod_{i\ge 1}
 \qbin{m_i^{(a)}+p_i^{(a)}}{m_i^{(a)}}
\end{equation}
where the sum is over all $\{m_i^{(a)}\in\N\mid 1\le a\le n,\; i\ge 1\}$
subject to the constraints
\begin{equation}\label{m restr}
\sum_{a=1}^n \sum_{i\ge 1} i m_i^{(a)} \alpha_a
 = \sum_{a=1}^n \sum_{i\ge 1} iL_i^{(a)} \La_a - \La.
\end{equation}
The variables $p_i^{(a)}$ and the exponent $\cc(\{m\})$ are defined as 
\begin{align}
\label{vac}
p_i^{(a)} &= \sum_{j\ge 1} L_j^{(a)} \min(i,j) - 
 \sum_{b=1}^n (\alpha_a|\alpha_b) \sum_{k\ge 1} \min(t_b i,t_a k) m_k^{(b)}\\
\label{cc}
\cc(\{m\}) &= \frac{1}{2} \sum_{a,b=1}^n (\alpha_a|\alpha_b)
 \sum_{j,k\ge 1} \min(t_b j,t_a k) m_j^{(a)} m_k^{(b)}
\end{align}
where $t_a=\frac{2}{(\alpha_a\mid\alpha_a)}$. Recall that
$(\cdot\mid\cdot)$ is normalized such that $(\alpha_a\mid\alpha_a)=2$
if $\alpha_a$ is a long root.
Then it was conjectured \cite{HKOTY} that
\begin{equation*}
\Xbar(B,\La)=\Fbar(B,\La).
\end{equation*}
For type $A_n^{(1)}$ and general $B$ this is proven in \cite{KSS} 
and for type $C_n^{(1)}$ and $B=\bigotimes_{a=1}^n 
(B^{a,1})^{\otimes L_1^{(a)}}$ a proof is given in \cite{OSS}.
Parts of the proofs given in \cite{KSS,OSS} are quite technical,
but we would like to highlight the general ideas of the proof 
which also give more insight into the fermionic formulas.

\subsection{Rigged configurations}

Rigged configurations provide a combinatorial interpretation
of the fermionic formula \eqref{fermi kostka}. We will focus
first on type $A_n^{(1)}$.

The sum over the variables $\{m_i^{(a)}\}$ in \eqref{fermi kostka}
subject to the restriction \eqref{m restr} can be interpreted as follows.
Let $\nu=(\nu^{(1)},\ldots,\nu^{(n)})$ be a sequence of partitions
with constraints on their sizes given by
\begin{equation}\label{nu restr}
|\nu^{(a)}|=-\sum_{j=1}^a \la_j+\sum_{i\ge 1}\sum_{b=1}^n i\;
 L_i^{(b)} \min(a,b)
\end{equation}
where $\la$ is the partition corresponding to the weight $\La$.
Then \eqref{nu restr} and \eqref{m restr} are equivalent provided
that $m_i^{(a)}$ is interpreted as the number of parts of $\nu^{(a)}$
of size $i$. In terms of $\nu$ the definitions \eqref{vac} and \eqref{cc}
read
\begin{align*}
P_i^{(a)}(\nu)&=Q_i(\nu^{(a-1)})-2Q_i(\nu^{(a)})+Q_i(\nu^{(a+1)})
 +\sum_{j\ge 1} L_j^{(a)} Q_i(j)\\
cc(\nu)&=\sum_{a=1}^n \sum_{i\ge 1} \alpha_i^{(a)}(\alpha_i^{(a)}
 -\alpha_i^{(a+1)})
\end{align*}
where $Q_i(\mu)$ is the number of boxes in the first $i$ columns of the
partition $\mu$, $\alpha_i^{(a)}$ is the size of the $i$-th column
of $\nu^{(a)}$ and $p_i^{(a)}=P_i^{(a)}(\nu)$.

To interpret \eqref{fermi kostka} combinatorially one uses the fact
the $q$-binomial coefficient $\qbins{m+p}{m}$ is the generating
function of partitions in a box of size $m\times p$. More precisely,
these are the partitions with at most $m$ parts each not exceeding
$p$. Hence \eqref{fermi kostka} can be restated as
\begin{equation*}
\Fbar(B,\La)=\sum_{(\nu,J)\in\RC(B,\La)} q^{\cc(\nu,J)}
\end{equation*}
where $\RC(B,\La)$ is the set of all $(\nu,J)$ where
$\nu$ is a sequence of partitions satisfying \eqref{nu restr}
and $J$ is a double sequence of partitions
\begin{equation*} 
J=\{J^{(a,i)}\}_{\substack{1\le a\le n\\ i\ge 1}}.
\end{equation*}
The partition $J^{(a,i)}$ has to fit in a box of size
$m_i^{(a)}(\nu) \times P_i^{(a)}(\nu)$ where $m_i^{(a)}(\nu)=m_i^{(a)}$
is the number of parts of size $i$ in $\nu^{(a)}$.
In particular this requires that $P_i^{(a)}(\nu)\ge 0$ for all
$i\ge 1$ and $1\le a\le n$. The exponent is defined as
\begin{equation*} 
\cc(\nu,J)=\cc(\nu)+\sum_{a=1}^n\sum_{i\ge 1} |J^{(a,i)}|.
\end{equation*}
The elements in $\RC(B,\La)$ are called rigged configurations.

Originally, rigged configurations were introduced in papers by Kerov, 
Kirillov and Reshetikhin \cite{KKR,KR} in their study of the XXX model using
Bethe Ansatz techniques. Rigged configurations index the solutions
to the Bethe Ansatz equations.

In \cite{KSS} the fermionic formula was proven by showing that
there is a statistic preserving bijection between classically
restricted paths and rigged configurations. 
It should be noted that for type $A_n^{(1)}$ the intrinsic energy
is equal to the energy, $\Dbar_B=\Ebar_B$.
\begin{theorem}\cite{KSS}\label{thm:kss}
For $B=\bigotimes_{a=1}^n \bigotimes_{i\ge 1} (B^{a,i})^{\otimes L_i^{(a)}}$
a crystal of type $A_n^{(1)}$, there is a bijection 
$\phi:\P'(B,\La)\to \RC(B,\La)$ such that for $b\in \P'(B,\La)$ we have
$\Ebar_B(b)=\cc(\theta\circ\phi(b))$.
Here $\theta:\RC(B,\La)\to\RC(B,\La)$ maps $(\nu,J)$ to
$(\nu,\tilde{J})$ where $\tilde{J}$ is obtained from $J$
by complementing each $J^{(a,i)}$ in the box of dimensions
$m_i^{(a)}(\nu)\times P_i^{(a)}(\nu)$.
\end{theorem}
The bijection $\phi$ is given explicitly in \cite{KSS}, 
\cite[Section 5.4]{OSS}.

The fermionic formula for type $C_n^{(1)}$ crystals of the form
$B_C=\bigotimes_{a=1}^n (B_C^{a,1})^{\otimes L_1^{(a)}}$ was proven in
\cite{OSS} using an embedding of type $C_n^{(1)}$ crystals
into type $A_{2n-1}^{(1)}$ crystals. Let
\begin{equation*}
\Psi(B_C^{r,1}) = \begin{cases} 
 B_A^{2n-r,1} \otimes B_A^{r,1} & \text{if $1\le r<n$}\\
 B_A^{n,2} & \text{if $r=n$.} \end{cases}
\end{equation*}
Baker \cite{B} showed that there is an embedding
$\Psi^{r,1}:B_C^{r,1}\to \Psi(B_C^{r,1})$.
It can be defined by requiring that the $U_q(C_n)$-highest
weight vector $u_C^{r,1}$ in $B_C^{r,1}$ is mapped to
$u_A^{2n-r,1}\otimes u_A^{r,1}$ where $u_A^{r,s}$ is the
$U_q(A_{2n-1})$-highest weight vector in $B_A^{r,s}$, and
\begin{align*}
\Psi^{r,1} \circ f_i^C &= f_{2n-i}^A \circ f_i^A \circ \Psi^{r,1}\\
\Psi^{r,1} \circ e_i^C &= e_{2n-i}^A \circ e_i^A \circ \Psi^{r,1}.
\end{align*}
For a tensor product, define
$\Psi_L:B_C^{r_L,1}\otimes\cdots\otimes B_C^{r_1,1}\to 
\Psi(B_C^{r_L,1})\otimes\cdots\otimes\Psi(B_C^{r_1,1})$
by $\Psi_L=\Psi^{r_L,1}\otimes\cdots\otimes\Psi^{r_1,1}$.

\begin{theorem}\cite{OSS}\label{thm:emb}
Let $B_C=B_C^{r_L,1}\otimes\cdots\otimes B_C^{r_1,1}$. The image
$\Image(\phi\circ\Psi_L)$ of $\phi\circ\Psi_L:
\P'(B_C,\cdot) \to \RC(\Psi(B_C),\cdot)$ is characterized by the
set of rigged configurations $(\nu,J)$ satisfying:
\begin{enumerate}
\item $(\nu,J)^{(k)}=(\nu,J)^{(2n-k)}$.
\item All parts of $\nu^{(n)}$ are even.
\item All riggings in $(\nu,J)^{(n)}$ are even.
\end{enumerate}
\end{theorem}

This characterization of the image of $\phi\circ\Psi_L$ suggests
the following definition of type $C$ rigged configurations.
Let $\la$ be a partition and 
$B_C=\bigotimes_{a=1}^n (B_C^{a,1})^{\otimes L_1^{(a)}}$,
and let $\nu=(\nu^{(1)},\ldots,\nu^{(n)})$ be a sequence of partitions
with the properties
\begin{equation}\label{nu C}
\begin{split}
&|\nu^{(a)}| = -\sum_{j=1}^a \la_j + \sum_{b=1}^n L_1^{(b)}
\min(a,b) \qquad \text{for $1\le a\le n$}\\
&\text{$\nu^{(n)}$ has only even parts.}
\end{split}
\end{equation}
Define the vacancy numbers as
\begin{equation}\label{P C}
\begin{split}
 P_i^{(a)}(\nu) &= Q_i(\nu^{(a-1)}) - 2Q_i(\nu^{(a)}) + Q_i(\nu^{(a+1)})
  + L_1^{(a)}  \quad \text{for $1\le a<n$,} \\
 P_i^{(n)}(\nu) &= Q_i(\nu^{(n-1)}) - Q_i(\nu^{(n)}) +
  \frac{1}{2} L_1^{(n)} Q_i(2).
\end{split}
\end{equation}

The set of rigged configurations of type $C$ corresponding to
a weight $\Lambda$ with associated partition $\la$ and crystal $B_C$, 
denoted by $\RC_C(B_C,\La)$, is given by $(\nu,J)$ where $\nu$ is a sequence
of partitions satisfying \eqref{nu C} and $J$ is a double
sequence of partitions 
\begin{equation*}
J=\{J^{(a,i)}\}_{\substack{1\le a\le n\\ i\ge 1}}
\end{equation*}
where $J^{(a,i)}$ is a partition in a box of size $m_i^{(a)}(\nu)
\times P_i^{(a)}(\nu)$ with $P_i^{(a)}(\nu)$ as in \eqref{P C}
and $m_i^{(a)}(\nu)$ the number of parts of $\nu^{(a)}$.

It is shown in \cite{OSS} that $\Dbar_A\circ \Psi_L=2\Dbar_C$.
Hence using Theorem \ref{thm:emb} the statistics of type $C$ rigged 
configurations becomes
\begin{equation*}
\begin{split}
\cc_C(\nu,J) &= \cc_C(\nu)+\sum_{a=1}^n\sum_{i\ge 1} |J^{(a,i)}| \\
\text{where}\quad \cc_C(\nu) &= \sum_{i\ge 1} \Bigl(
\sum_{a=1}^{n-1} \alpha_i^{(a)}(\alpha_i^{(a)}-\alpha_i^{(a+1)})
+\frac{1}{2} \alpha_i^{(n) \;2} \Bigr)
\end{split}
\end{equation*}
which implies that
\begin{equation*}
\Xbar(B_C,\Lambda) = \sum_{(\nu,J)\in\RC_C(B_C,\Lambda)} q^{\cc_C(\nu,J)}.
\end{equation*}
It is also not so hard to show that
\begin{equation*}
\Fbar(B_C,\Lambda) = \sum_{(\nu,J)\in\RC_C(B_C,\Lambda)} q^{\cc_C(\nu,J)}
\end{equation*}
by identifying
\begin{equation*}
\begin{split}
p_i^{(a)}&=\begin{cases}
 P_i^{(a)}(\nu) & \text{for $1\le a<n$}\\
 P_{2i}^{(n)}(\nu) & \text{for $a=n$} \end{cases}\\
m_i^{(a)}&=\begin{cases}
 m_i(\nu^{(a)}) & \text{for $1\le a<n$}\\
 m_{2i}(\nu^{(n)}) & \text{for $a=n$.} \end{cases}
\end{split}
\end{equation*}
This proves that $\Xbar(B_C,\Lambda)=\Fbar(B_C,\Lambda)$.

\subsection{Level-restricted case}

Fermionic formulas for the level-restricted one-dimensional
configuration sums $\Xbar^\ell(B,\Lambda)$ were conjectured in \cite{HKOTY}
for all $\gggh$ as in section \ref{sec:classic} and
special weight $\Lambda=0$.
Let $\ell\in\N$ and define the following polynomial
in $q$ depending on $B$
\begin{equation}\label{fermi lev kostka}
\Fbar^\ell(B)=\sum_{\{m\}} q^{\cc^\ell(\{m\})} \prod_{(a,i)\in H^\ell}
 \qbin{m_i^{(a)}+p_i^{(a)}}{m_i^{(a)}}
\end{equation}
where $H^\ell=\{(a,i)\mid 1\le a\le n, 1\le j\le t_a\ell\}$ and the sum 
is over all $\{m_i^{(a)}\in\N\mid (a,i)\in H^\ell\}$ subject to the 
constraints
\begin{equation}\label{m lev restr}
\sum_{(a,i)\in H^\ell} i\; m_i^{(a)} \alpha_a
 = \sum_{(a,i)\in H^\ell} i\; L_i^{(a)} \La_a.
\end{equation}
The variables $p_i^{(a)}$ and the exponent $\cc^\ell(\{m\})$ are defined as 
\begin{align}
\label{p ell}
p_i^{(a)} &= \sum_{j=1}^{t_a\ell} L_j^{(a)} \min(i,j) - 
 \sum_{(b,k)\in H^\ell} (\alpha_a|\alpha_b) \min(t_b i,t_a k) m_k^{(b)}\\
\label{cc ell}
\cc^\ell(\{m\}) &= \frac{1}{2} \sum_{(a,j),(b,k)\in H^\ell} 
 (\alpha_a|\alpha_b) \min(t_b j,t_a k) m_j^{(a)} m_k^{(b)}.
\end{align}
Then it was conjectured \cite{HKOTY} that
\begin{equation}\label{conj level}
\Xbar^\ell(B,0)=\Fbar^\ell(B).
\end{equation}

\subsection{Level-restricted case: type $A_n^{(1)}$}

For type $A_n^{(1)}$ the conjecture \eqref{conj level}
and its generalization to arbitrary weights $\Lambda$ was proven 
in \cite{SSa}. These formulas can again be understood in terms
of rigged configurations. We will explain this here since it will
enable us to derive the level-restricted fermionic formulas
for type $C_n^{(1)}$ in the next section.

Let $\la=(\la_1,\ldots,\la_{n+1})$ be the partition corresponding
to the dominant integral weight $\Lambda$.
A partition $\la$ is restricted of level $\ell$
if $\la_1-\la_{n+1}\le \ell$. Define $\lt = \ell-(\la_1-\la_{n+1})$, 
which is nonnegative by assumption.
Set $\la'=(\la_1-\la_{n+1},\ldots,\la_n-\la_{n+1})^t$ (where ${}^t$ stands
for transpose) and denote the set of all column-strict tableaux of shape 
$\la'$ over the alphabet $\{1,2,\ldots,\la_1-\la_{n+1}\}$ by $\CST(\la')$. 
Define a table of modified vacancy numbers depending on a sequence of
partitions $\nu$ and $t\in\CST(\la')$ by
\begin{equation} \label{t vacancy}
  P_i^{(a)}(\nu,t) =
  P_i^{(a)}(\nu) - \sum_{j=1}^{\la_a-\la_{n+1}} \chi(i\ge\lt+t_{j,a})
+ \sum_{j=1}^{\la_{a+1}-\la_{n+1}} \chi(i\ge\lt+t_{j,a+1})
\end{equation}
where $t_{j,a}$ denotes the entry in $t$ in row $j$ and column $a$
and $\chi(S)=1$ if the statement $S$ is true and $\chi(S)=0$ otherwise.

\begin{definition}\label{def:levrc}
Let $B=\bigotimes_{a=1}^n \bigotimes_{i=1}^\ell (B^{a,i})^{\otimes L_i^{(a)}}$
be a crystal of type $A_n^{(1)}$ and $\La$ a dominant integral weight
with corresponding partition $\la$.
Say that $(\nu,J)\in\RC(B,\La)$ is restricted of level $\ell$
provided that
\begin{enumerate}
\item $\nu_1^{(a)} \le \ell$ for all $a$.
\item There exists a tableau $t\in\CST(\la')$, such that for every 
$i,a\ge 1$, the largest part of $J^{(a,i)}$ does not exceed 
$P_i^{(a)}(\nu,t)$.
\end{enumerate}
Denote by $\RC^\ell(B,\La)$ the set of $(\nu,J)\in\RC(B,\La)$ that are 
restricted of level $\ell$.
\end{definition}
Note in particular that the second condition requires that
$P_i^{(a)}(\nu,t)\ge 0$ for all $i,a\ge 1$.

\begin{example}\label{ex:vacuum}
Let $\la=(m^{n+1})$ be rectangular with $n+1$ rows.
Then $\la'=\emptyset$ and $P_i^{(a)}(\nu,\emptyset)=P_i^{(a)}(\nu)$
for all $i,a\ge 1$ so that the modified vacancy numbers are equal
to the vacancy numbers.
\end{example}

\begin{theorem}\cite[Theorem 8.2]{SSa}
\label{thm:lev image A}
The bijection $\phi:\P'(B,\La)\to \RC(B,\La)$ restricts
to a well-defined bijection $\phi:\P^\ell(B,\La)\to \RC^\ell(B,\La)$.
\end{theorem}

Since $\Ebar_B=\cc\circ\theta\circ\phi$ by Theorem \ref{thm:kss} it 
follows from Theorem \ref{thm:lev image A} that for type $A_n^{(1)}$
\begin{equation}\label{lev comb A}
\Xbar^\ell(B,\La)=\sum_{(\nu,J)\in \RC^\ell(B,\La)} q^{\cc(\theta(\nu,J))}.
\end{equation}

Note that $\cc(\theta(\nu,J))=\cc(\nu)+\sum_{a=1}^n \sum_{i=1}^\ell
 P_i^{(a)}(\nu)m_i^{(a)}(\nu)-\sum_{a=1}^n \sum_{i=1}^\ell |J^{(a,i)}|$.
Let $\SCST(\la')$ be the set of all nonempty subsets of $\CST(\la')$.
Since the $q$-binomial $\qbins{m+p}{m}$ is the generating function of 
partitions with at most $m$ parts each not exceeding $p$
it follows by inclusion-exclusion that
\begin{equation}\label{lev type A}
\Xbar^\ell(B,\La)=\sum_{S\in\SCST(\la')} (-1)^{|S|+1}
 \sum_{\{m\}} q^{\ch^\ell(\{m\})} \prod_{(a,i)\in H^\ell}
 \qbin{m_i^{(a)}+p_i^{(a)}(S)}{m_i^{(a)}}_{1/q}
\end{equation}
where the sum is over all $\{m_i^{(a)}\in\N\mid (a,i)\in H^\ell\}$ subject 
to the constraints
\begin{equation*}
\sum_{(a,i)\in H^\ell} i m_i^{(a)} \alpha_a
 = \sum_{(a,i)\in H^\ell} i L_i^{(a)} \La_a-\La.
\end{equation*}
Also $p_i^{(a)}(S)$ and $\ch^\ell(\{m\})$ are defined as
\begin{align*}
p_i^{(a)}(S)&=p_i^{(a)}+\min_{t\in S}\Bigl\{
 - \sum_{j=1}^{\la_a-\la_{n+1}} \chi(i\ge\lt+t_{j,a})
 + \sum_{j=1}^{\la_{a+1}-\la_{n+1}} \chi(i\ge\lt+t_{j,a+1})\Bigr\}\\
\ch^\ell(\{m\})&=\cc^\ell(\{m\})+\sum_{(a,i)\in H^\ell} p_i^{(a)} m_i^{(a)}
\end{align*}
with $p_i^{(a)}$ as in \eqref{p ell} and $\cc^\ell(\{m\})$ as in
\eqref{cc ell}. $\qbin{m+p}{m}_{1/q}$ is the $q$-binomial with $q$
replaced by $1/q$. In particular if $\la=(m^{n+1})$ as in
Example \ref{ex:vacuum} so that $\Lambda=0$ and $p_i^{(a)}(\nu,\{\emptyset\})
=p_i^{(a)}$ the fermionic form \eqref{lev type A} reduces to 
\eqref{fermi lev kostka} since $\qbin{m+p}{m}_{1/q}=q^{-mp}\qbin{m+p}{m}$. 
Further details can be found in \cite{SSa}.

It should be remarked that even though \eqref{lev type A} contains
explicit signs, it is clear from the equivalent combinatorial 
formula \eqref{lev comb A} that it is a nonnegative polynomial in $q$.

\subsection{Level-restricted case: type $C_n^{(1)}$}
\label{sec:lev C}

In this section we will show how the level-restricted
fermionic formulas for type $C_n^{(1)}$ can be obtained from
Theorems \ref{thm:lev image A} and \ref{thm:emb}.

Under the embedding $\Psi$ a dominant weight 
$\Lambda_C=\sum_{k=1}^m \Lambda_{i_k}^C$ of type $C_n^{(1)}$ becomes the 
weight $\Lambda_A=\sum_{k=1}^m (\Lambda_{i_k}^A+\Lambda_{2n-i_k}^A)$
of type $A_{2n-1}^{(1)}$ where all $1\le i_k\le n$. In terms of the 
corresponding partitions $\la^A$ and $\la^C$ this implies that
\begin{align*}
\la_a^A-\la_{2n}^A&=\la_1^C+\la_a^C && \text{for $1\le a\le n$}\\
\la_a^A-\la_{2n}^A&=\la_1^C-\la_{2n+1-a}^C && \text{for $n< a\le 2n$.}
\end{align*}
Let $B_C=\bigotimes_{a=1}^n (B_C^{a,1})^{\otimes L_1^{(a)}}$.
Hence, under the embedding $\Psi_L:B_C\to \Psi(B_C)$, the conditions
for level-restriction for rigged configurations of type $A$ as given
in Definition \ref{def:levrc} become the following.

For a partition $\la^C$ define $(\la^C)'=(2\la_1^C,\la_1^C+\la_2^C,
\ldots,\la_1^C+\la_n^C,\la_1^C-\la_n^C,\la_1^C-\la_{n-1}^C,\ldots,
\la_1^C-\la_1^C)^t$. Let $\CST((\la^C)')$ be the set of all semi-standard
tableaux of shape $(\la^C)'$ in the alphabet $\{1,2,\ldots,2\la_1^C\}$.
For $t\in \CST((\la^C)')$ set
\begin{equation*}
f_i^{(a)}(t)= \begin{cases}
 -\sum_{j=1}^{\la_1^C+\la_a^C} \chi(i\ge 2\ell-2\la_1^C+t_{j,a})&\\
 \quad+\sum_{j=1}^{\la_1^C+\la_{a+1}^C} \chi(i\ge 2\ell-2\la_1^C+t_{j,a+1})
 &\text{for $1\le a<n$}\\
 -\sum_{j=1}^{\la_1^C+\la_n^C} \chi(i\ge 2\ell-2\la_1^C+t_{j,n})&\\
 \quad+\sum_{j=1}^{\la_1^C-\la_n^C} \chi(i\ge 2\ell-2\la_1^C+t_{j,n+1})
 &\text{for $a=n$}\\
 -\sum_{j=1}^{\la_1^C-\la_{2n+1-a}^C} \chi(i\ge 2\ell-2\la_1^C+t_{j,a})&\\
 \quad+\sum_{j=1}^{\la_1^C-\la_{2n-a}^C} \chi(i\ge 2\ell-2\la_1^C+t_{j,a+1})
 &\text{for $n<a<2n$.}
\end{cases}
\end{equation*}
Define modified vacancy numbers as
\begin{equation}\label{P C lev}
P_i^{(a)}(\nu,t)=P_i^{(a)}(\nu)+\begin{cases}
 \min\{f_i^{(a)}(t),f_i^{(2n-a)}(t)\} & \text{for $1\le a<n$}\\
 \frac{1}{2}f_i^{(n)}(t) & \text{for $a=n$}
\end{cases}
\end{equation}
with $P_i^{(a)}(\nu)$ as defined in \eqref{P C}.

\begin{definition}
Let $\La_C$ be a dominant weight and $\la^C$ the corresponding partition.
Let $B_C=\bigotimes (B_C^{a,1})^{\otimes L_1^{(a)}}$
be a crystal of type $C_n^{(1)}$.
Say that $(\nu,J)\in\RC_C(B_C,\La_C)$ is restricted of level $\ell$
provided that
\begin{enumerate}
\item $\nu_1^{(a)} \le 2\ell$ for all $a$.
\item There exists a tableau $t\in\CST((\la^C)')$, such that for every 
$i,a\ge 1$, the largest part of $J^{(a,i)}$ does not exceed 
$P_i^{(a)}(\nu,t)$ defined in \eqref{P C lev}.
\end{enumerate}
Denote by $\RC_C^\ell(B_C,\La_C)$ the set of $(\nu,J)\in\RC_C(B_C,\La_C)$ 
that are restricted of level $\ell$.
\end{definition}

It follows that
\begin{equation*}
\Xbar^\ell(B_C,\La_C)=\sum_{(\nu,J)\in \RC_C^\ell(B_C,\La_C)} 
 q^{\cc_C(\theta(\nu,J))}.
\end{equation*}

Let $\SCST((\la^C)')$ be the set of all nonempty subsets of $\CST((\la^C)')$.
By the same arguments as in the type $A$ we find
\begin{equation}\label{lev type C}
\Xbar^\ell(B_C,\La_C)=\sum_{S\in\SCST((\la^C)')} (-1)^{|S|+1}
 \sum_{\{m\}} q^{\ch^\ell(\{m\})} \prod_{(a,i)\in H^\ell}
 \qbin{m_i^{(a)}+p_i^{(a)}(S)}{m_i^{(a)}}_{1/q}
\end{equation}
where the sum is over all $\{m_i^{(a)}\in\N\mid (a,i)\in H^\ell\}$ subject 
to the constraints
\begin{equation*}
\sum_{(a,i)\in H^\ell} i m_i^{(a)} \alpha_a
 = \sum_{(a,1)\in H^\ell} L_1^{(a)} \La_a^C-\La_C.
\end{equation*}
The variable $p_i^{(a)}(S)$ is defined as
\begin{align*}
p_i^{(a)}(S)&=p_i^{(a)}+\begin{cases}
 \min_{t\in S}\{f_i^{(a)}(t),f_i^{(2n-a)}(t)\} & \text{for $1\le a<n$}\\
 \frac{1}{2}\min_{t\in S}\{f_i^{(n)}(t)\} & \text{for $a=n$}
\end{cases}\\
\ch^\ell(\{m\})&=\cc^\ell(\{m\})+\sum_{(a,i)\in H^\ell} p_i^{(a)} m_i^{(a)}
\end{align*}
with $p_i^{(a)}$ as in \eqref{vac} and $\cc^\ell(\{m\})$ as in
\eqref{cc}. In particular if $\la^C=\emptyset$ 
so that $\Lambda_C=0$ the fermionic form \eqref{lev type C} reduces 
to \eqref{fermi lev kostka}.

\section{Summary and open problems}
\label{sec:open}

Equating the bosonic and fermionic evaluations for
$\Xbar(B,\La)$ and $\Xbar^\ell(B,\La)$ as given in sections
\ref{sec:bose} and \ref{sec:fermi} yields polynomial identities in $q$.
One may take limits of these identities in various ways to obtain
$q$-series identities. We refer the interested reader to
\cite{HKKOTY,HKOTY,SSa,SW} for details.

It is still an outstanding problem to prove the conjectured fermionic
formulas \eqref{fermi kostka} and \eqref{fermi lev kostka} for general 
$\gggh$. There is strong evidence that the rigged configuration approach
will still be applicable in this case. Part of this program
also requires the proof of the existence and structure of the crystals
$B^{r,s}$ for general $\gggh$ as mentioned in Conjecture \ref{conj:Brs}.

The astute reader will have noticed that the energy
function introduced in section \ref{sec:crystal}
does not actually specialize to the 
statistics on path of section \ref{sec:HHM} which yield 
the Rogers--Ramanujan identities. This is due to the fact that
the Rogers--Ramanujan identities correspond to a
non-unitary physical model (more precisely the 
minimal model $M(2,5)$) whereas the crystal base
theory introduced in section \ref{sec:crystal}
is related to unitary physical models. To really
embed the Rogers--Ramanujan identities into the framework
of crystal base theory it is hence necessary to generalize the
definition of the (intrinsic) energy function of section \ref{sec:crystal} 
to include the non-unitary setting. For $\gggh=\widehat{\mathfrak{sl}}_2$ 
results are available \cite{FB,BMS,FW,FLPW}. However for general $\gggh$ 
details are not yet known in general. Some results in this direction can 
be found in \cite{N}.

\bibliographystyle{amsplain}

\end{document}